\documentclass[final,1p]{elsarticle}

\usepackage{amssymb,latexsym}
\usepackage{amsmath}

\journal{: \quad  Int. J. Math. Educ. Sci. Technol.}

\begin{document}

\begin{frontmatter}

\title{A shortcut for evaluating some log integrals from products and limits}

\author{F. M. S. Lima}

\address{Institute of Physics, University of Brasilia, P.O. Box 04455, 70919-970, Brasilia-DF, Brazil}
%{\rm Present address:} Departamento de Fisica, Universidade Federal de Pernambuco, 50670-901, Recife-PE, Brazil}

%Tel.: +55 (061) 9622-5871  --  FAX: +55 (061) 3307-2363

\ead{fabio@fis.unb.br}

%\date{\today}

\begin{abstract}
In this short paper, I introduce an elementary method for exactly evaluating the definite integrals $\, \int_0^{\pi}{\ln{(\sin{\theta})}\,d\theta}$, $\int_0^{\pi/2}{\ln{(\sin{\theta})}\,d\theta}$, $\int_0^{\pi/2}{\ln{(\cos{\theta})}\,d\theta}$, and $\int_0^{\pi/2}{\ln{(\tan{\theta})}\,d\theta} \,$ in finite terms.  The method consists in to manipulate the sums obtained from the logarithm of certain products of trigonometric functions at rational multiples of $\pi$, putting them in the form of Riemann sums. As this method does not involve any search for primitives, it represents a good alternative to more involved integration techniques.  As a bonus, I show how to apply the method for easily evaluating $\,\int_0^1{\ln{\Gamma(x)} \, d x}$.
\newline
\end{abstract}


\end{frontmatter}

\section{Introduction}
The simplicity of the functions $\,\ln{(\sin{\theta})}$, $\ln{(\cos{\theta})}$, and $\ln{(\tan{\theta})}$, as well as the fact they are continuous (even differentiable) except at some isolated points, suggests that the evaluation of the definite integrals $\int_0^{\pi}{\ln{(\sin{\theta})}\,d\theta}$, $\int_0^{\pi/2}{\ln{(\sin{\theta})}\,d\theta}$, $\int_0^{\pi/2}{\ln{(\cos{\theta})}\,d\theta}$, and $\int_0^{\pi/2}{\ln{(\tan{\theta})}\,d\theta}$  via the Fundamental Theorem of Calculus should be a straightforward task,\footnote{~These integrals are, in fact, \emph{improper} because each integrand has at least one infinite discontinuity in the integration interval.} which is not true.\footnote{~Note that it is not possible to express the corresponding indefinite integrals in a finite, closed-form expression involving only \emph{elementary functions}.}  The usual methods for the analytic evaluation of such integrals involve advanced integration techniques, such as to expand the integrand in a power series and then integrate it term-by-term (see Sec.~11.9 of Ref.~\cite{Stewart}), or the evaluation of a suitable contour integral on the complex plane in view to apply the Cauchy's residue theorem (see Secs.~4.1 and 4.2 of Ref.~\cite{ComplexBook}).  However, these methods present some disadvantages when applied to definite integrals of ``log-trig'' functions.  The expansion of the integrand in a Maclaurin series, e.g., is not possible for $\ln{(\sin{\theta})}$ and $\ln{(\tan{\theta})}$. Though this series expansion is possible for $\ln{(\cos{\theta})}$, it is difficult to determine a closed-form expression for the general term and then to recognize the number it represents.  The evaluation of a contour integral on the complex plane has the inconveniences of requiring the choice of a \emph{suitable} integration path, usually a difficult task, and yielding logarithms of complex (non-real) numbers and/or non-elementary transcendental functions (e.g., dilogarithm, elliptic, and hypergeometric functions),  which often makes it obscure the final result since these functions are either unknown or badly-known for most students. These inconveniences are just what one finds by appealing to mathematical softwares. For instance, Maple (release 13) and Mathematica (release 7) both return the following `stodgy' result for $\int{\ln{(\sin{x})}\,dx}$: %, with $0<x<\pi$:
\begin{eqnarray*}
\frac{i}{2} \, {\rm Li}_2 \left( { e^{2\,ix}} \right) +i\,\frac{x^2}{2} +x\, \ln{\left(\sin{x}\right)} -x\,\ln{\left( 1 -e^{2\,ix} \right)} ,
\end{eqnarray*}
where $\mathrm{Li}_2(x) := \sum_{n=1}^{\infty}{\dfrac{x^n}{n^2}}$ is the dilogarithm function.

In this work, the integrals $\int_0^{\pi}{\ln{(\sin{\theta})}\,d\theta}$, $\int_0^{\pi/2}{\ln{(\sin{\theta})}\,d\theta}$, $\int_0^{\pi/2}{\ln{(\cos{\theta})}\,d\theta}$, and $\int_0^{\pi/2}{\ln{(\tan{\theta})}\,d\theta}$ are easily evaluated by applying the natural logarithm to certain products of trigonometric functions at rational multiples of $\pi$, which yields sums that can be written in the form of Riemann sums. The closed-form expressions emerge when we take suitable limits (as the number of terms tends to infinity), without any search for primitives.  I also show how the method can be used for easily evaluating $\,\int_0^1{\ln{\Gamma(x)} \, d x}$.  %, together with the reflection property of the gamma function,

\section{Some products of trigonometric functions}
Let me present the products of trigonometric functions at rational multiples of $\pi$ that will serve as the basis for my method.%, as described in the next section.

In Appendix A.3 of Ref.~\cite{Niven}, in presenting an elementary proof for the Euler result $\,\sum_{n=1}^{\infty}{{\,1/n^2}} = {\,\pi^2/6}$, the authors prove an identity involving $\cot^2{\theta}$, $\theta$ being a rational multiple of $\pi$.
%Their proof follows from a comparison of the De Moivre's theorem --- namely, $(\cos{\theta} + i \, \sin{\theta})^m = \cos{m \theta} + i \sin{m \theta}$, valid for all positive integers $m$ --- and the binomial theorem for $\cos{m \theta} + i \sin{m \theta} = \sin^m{\theta} \, \left(\cot{\theta}+i\right)^m$. They put $m=2 N+1$ and show that
%\begin{eqnarray*}
%\sin{(2 N+1)\theta} = \sin^{2 N+1}{\theta} \: \times \: F(\cot^2{\theta}) \, ,
%\end{eqnarray*}
%where
%Note that, for all $\theta_n = {\,n \pi/(2N+1)}$, $n=1,\ldots,N$, one has $\sin{(2N+1)\theta_n} = 0$, but $\sin{\theta_n} \ne 0$. Then, the roots of the polynomial equation $F(x)=0$ are just the numbers $\cot^2{\theta_n}$, $n=1,\ldots,N$.
By applying the well-known rule for the product of roots of a polynomial equation $\,F(x) = 0\,$ to
\begin{eqnarray*}
F(x) = \sum_{j=0}^N{(-1)^j\,\binom{2N+1}{2j+1}\,x^{N-j}} = (2N+1)\,\prod_{n=1}^N{\left(x-\cot^2{\theta_n}\right)} \, ,
%\binom{2 N+1}{1}x^N - \binom{2 N+1}{3}x^{N-1} + \binom{2 N+1}{5}x^{N-2} - \ldots + (-1)^N .
\end{eqnarray*}
where $\theta_n = {\,n \pi/(2N+1)}$, $n=1,\ldots,N$, it is shown that
\begin{eqnarray*}
\prod_{n=1}^{N}{\cot^2{\theta_n}} = (-1)^N \, \frac{(-1)^N}{\binom{2 N+1}{1}} = \frac{1}{2 N+1} \, .
\end{eqnarray*}
By extracting the square-root of the inverse of each side, one finds
\begin{equation}
\prod_{n=1}^{N}{\tan{\left(\frac{n \, \pi}{2 N+1}\right)}} = \sqrt{2N+1} \, .
\label{eq:tan}
\end{equation}

In that appendix, we also find a proof for the following identity:
\begin{equation}
\prod_{n=1}^{N-1} \sin{\left(\frac{\pi n}{N}\right)} = \frac{N}{2^{N-1}} \, .
\label{eq:seno}
\end{equation}

From the symmetry relation $\sin{(\pi/2 + \alpha)} = \sin{(\pi/2 - \alpha)}$, it is easy to deduce that
\begin{eqnarray*}
\prod_{n=1}^{N-1} \sin{\left(\frac{\pi n}{N}\right)} = \prod_{n=1}^{\left\lfloor\frac{N}{2}\right\rfloor} \sin{\left(\frac{\pi n}{N}\right)} \: \sin{\left(\frac{\pi (N-n)}{N}\right)} = \prod_{n=1}^{\left\lfloor\frac{N}{2}\right\rfloor} \sin{\left(\frac{\pi n}{N}\right)} \: \sin{\left(\pi -\frac{\pi n}{N}\right)} \, ,
\end{eqnarray*}
valid for all $N>1$, which implies that
\begin{equation}
\prod_{n=1}^{\left\lfloor\frac{N}{2}\right\rfloor} \sin^2{\left(\frac{\pi n}{N}\right)} = \frac{N}{2^{N-1}} \, .
\label{eq:halfseno}
\end{equation}

For a product of cosines, note that, $\forall \; \alpha \in \left[0, {\,\pi/2}\right]$, $\sin{\alpha} = \cos{({\,\pi/2}-\alpha)}$.  By taking $\alpha = {\,\pi n/N}$ and applying this trigonometric identity on Eq.~\eqref{eq:halfseno}, one readily finds
\begin{equation}
\prod_{n=1}^{\left\lfloor\frac{N}{2}\right\rfloor} \cos^2{\left(\frac{\pi}{2} -\frac{\pi n}{N}\right)} = \frac{N}{2^{N-1}} \, .
\label{eq:cos}
\end{equation}

These are the trigonometric products needed for evaluating the integrals we are interested in here.

\section{Evaluation of definite integrals from products}
The general idea underlying my method is to take the logarithm of a product of positive terms, convert it into a sum of logarithms and then to put this sum in the form of a Riemann sum with equally-spaced subintervals, whose limit as the number of terms tends to infinity is a definite integral.  Let us apply this procedure to the products of trigonometric functions in Eqs.~\eqref{eq:tan}--\eqref{eq:cos}.

For instance, by taking the logarithm of each side of Eq.~\eqref{eq:seno}, one has
\begin{equation}
\sum_{n=1}^{N-1}{\ln{\left[\sin{\left(\pi \, \frac{n}{N}\right)}\right]}} = \ln{N} - (N-1)\,\ln{2} \, .
\label{eq:lnseno}
\end{equation}
By dividing both sides by $N-1$, one gets
\begin{equation}
\sum_{n=1}^{N-1}{\frac{\ln{\left[\sin{\left(\pi \, \frac{n}{N}\right)}\right]}}{N-1}} = \frac{\ln{N}}{N-1} - \ln{2} \, .
\label{eq:aux1}
\end{equation}
Now, let us define $x_n = {\,n/N}$ and $\Delta x = {\,1/(N-1)}$. Equation~\eqref{eq:aux1} then reads
\begin{equation}
\sum_{n=1}^{N-1} {\ln{\left[\sin{\left(\pi \, x_n\right)}\right]}\,\Delta x} = \frac{\ln{N}}{N-1} - \ln{2} \, .
\label{eq:aux2}
\end{equation}
Clearly, the sum at the left-hand side has the form of a Riemann sum in which the grid points $x_n$ are equally spaced by $\Delta x$. By taking the limit as $N \rightarrow \infty$ on both sides of this equation and noting that $\lim_{N \rightarrow \infty}{\dfrac{\ln{N}}{N-1}} = 0$, which follows from L'Hopital rule, one has
\begin{equation}
\lim_{N \rightarrow \infty}{\sum_{n=1}^{N-1} {\ln{\left[\sin{\left(\pi \, x_n\right)}\right]}\,\Delta x}} = - \ln{2} \, ,
\label{eq:aux3}
\end{equation}
which means that
\begin{equation}
\int_{0^{+}}^{1^{-}}{\ln{\,\sin{\left(\pi \, x\right)}} \: d x} = - \ln{2} \, .
\label{eq:int1}
\end{equation}
The change of variable $\theta=\pi\,x$ promptly yields
\begin{equation}
\int_0^{\,\pi}{\ln{\,\sin{\theta}} \: d \theta} = - \,\pi \, \ln{2} \, .
\label{eq:int2}
\end{equation}
%This result can also be found by applying the trick substitutions found in Sec.~12.5 of Ref.~\cite{Irresistiveis}.
%\begin{flushright} $\Box$ \end{flushright}

When the above procedure is applied to Eq.~\eqref{eq:halfseno} one finds that
%\begin{equation}
%\sum_{n=1}^{\lfloor N/2 \rfloor}{2 \, \ln{\left[\sin{\left(\pi \, \frac{n}{N}\right)}\right]}} = \ln{N} - (N-1)\,\ln{2} \, .
%\label{eq:lnhalf}
%\end{equation}
%By dividing both sides by $N-1$ and substituting $x_n$ and $\Delta x$ as before, one gets
\begin{eqnarray*}
2 \sum_{n=1}^{\lfloor N/2 \rfloor} {\ln{\left[\sin{\left(\pi \, x_n\right)}\right]}\,\Delta x} = \frac{\ln{N}}{N-1} - \ln{2} \, .
%\label{eq:aux5}
\end{eqnarray*}
The limit as $N \rightarrow \infty$ yields
\begin{eqnarray*}
2 \, \int_0^{\frac12}{\ln{\,\sin{\left(\pi \, x\right)}} \: d x} = -\ln{2} \, .
%\label{eq:int3}
\end{eqnarray*}
The change of variable $\theta=\pi\,x$ yields
\begin{equation}
\int_0^{\,\pi/2}{\ln{\,\sin{\theta}} \: d \theta} = - \, \frac{\pi}{2} \, \ln{2} \, .
\label{eq:int4}
\end{equation}
%\begin{flushright} $\Box$ \end{flushright}

Now, from the product of cosines in Eq.~\eqref{eq:cos} one has
%\begin{eqnarray*}
%\sum_{n=1}^{\lfloor N/2 \rfloor}{2 \, \ln{\left[\cos{\left(\frac{\pi}{2} -\pi \, \frac{n}{N}\right)}\right]}} = \ln{N} - (N-1) \,\ln{2} \, .
%\label{eq:lncos}
%\end{eqnarray*}
%Again, by dividing both sides by $N-1$, one gets
\begin{eqnarray*}
2 \, \sum_{n=1}^{\lfloor N/2 \rfloor} {\ln{\left[\cos{\left(\frac{\pi}{2} -\pi \, x_n\right)}\right]}\,\Delta x} = \frac{\ln{N}}{N-1} - \ln{2} \, .
%\label{eq:aux7}
\end{eqnarray*}
The limit as $N \rightarrow \infty$ yields
\begin{eqnarray*}
2 \, \int_0^{\frac12}{\ln{\,\cos{\left(\frac{\pi}{2}-\pi \, x\right)}} \: d x} = -\ln{2} \, .
%\label{eq:int5}
\end{eqnarray*}
The change of variable $\theta=\frac{\pi}{2}-\pi\,x$ yields
\begin{equation}
\int_0^{\,\pi/2}{\ln{\,\cos{\theta}} \: d \theta} = - \, \frac{\pi}{2} \, \ln{2} \, .
\label{eq:int6}
\end{equation}
%\begin{flushright} $\Box$ \end{flushright}

From the product of tangents in Eq.~\eqref{eq:tan}, one has
\begin{equation*}
\sum_{n=1}^{N}{\ln{\left[\tan{\left(\pi \, \frac{n}{2 N+1}\right)}\right]}} = \frac{1}{2} \, \ln{(2 N+1)} \, .
%\label{eq:lntan}
\end{equation*}
By substituting $2 N+1 = M$ (hence $M$ is an odd positive integer) and then dividing both sides by $M$, one has
\begin{eqnarray*}
\sum_{n=1}^{(M-1)/2} {\ln{\left[\tan{\left(\pi \, x_n\right)}\right]}\,\Delta x} = \frac{\ln{M}}{2 M} \, ,
%\label{eq:aux8}
\end{eqnarray*}
where $x_n={\,n/M}$ and $\Delta x = {\,1/M}$. %Again, the sum at the left-hand side has the form of a Riemann sum with equally-spaced grid points.
The limit as $M \rightarrow \infty$ yields
\begin{eqnarray*}
\lim_{M \rightarrow \infty}{\sum_{n=1}^{(M-1)/2} {\ln{\left[\tan{\left(\pi \, x_n\right)}\right]}\,\Delta x}} = \lim_{M \rightarrow \infty}{\frac{\ln{M}}{2 M}} = \frac12 \, \lim_{M \rightarrow \infty}{\frac{1}{M}} = 0 \, ,
%\label{eq:aux9}
\end{eqnarray*}
which means that
\begin{eqnarray*}
\int_0^{\frac12}{\ln{\,\tan{\left(\pi \, x\right)}} \: d x} = 0 \, .
%\label{eq:int7}
\end{eqnarray*}
The change of variable $\theta=\pi\,x$ yields
\begin{equation}
\int_0^{\,\pi/2}{\ln{\,\tan{\theta}} \: d \theta} = 0 \, .
\label{eq:int8}
\end{equation}
%\begin{flushright} $\Box$ \end{flushright}

Surprisingly, this method also works in evaluating the  `impossible' integral $\int_0^1{\ln{\Gamma(x)} \, d x}$, where $\Gamma(x)$ is the Euler gamma function.\footnote{Here, the word \emph{impossible} certainly reflects the opinion of most students.}  For this, let us add $(N-1) \, \ln{\pi}$ on both sides of Eq.~\eqref{eq:lnseno}, which yields
\begin{eqnarray*}
(N-1) \, \ln{\pi} - \sum_{n=1}^{N-1}{\ln{\left[\sin{\left(\pi \, \frac{n}{N}\right)}\right]}} = (N-1) \, \ln{\pi} -\ln{N} + (N-1)\,\ln{2} \, .
%\label{eq:lnnew}
\end{eqnarray*}
This promptly simplifies to %$\sum_{n=1}^{N-1}{\left[\ln{\pi}-\ln{\sin{\left(\pi \, x_n\right)}}\right]} = (N-1) \, \ln{(2 \, \pi)} -\ln{N}$
\begin{equation}
\sum_{n=1}^{N-1}{\ln{\left[\frac{\pi}{\sin{\left(\pi \, x_n\right)}}\right]}} = (N-1) \, \ln{(2 \, \pi)} -\ln{N} \, .
\label{eq:boa}
\end{equation}
%where $x_n = {\,n/N}$.
Now, let us make use of the reflection property $\Gamma{(x)} \cdot \Gamma{(1-x)} = {\, \pi / \sin{(\pi \, x)}}$, valid for all $x \not \in \mathbb{Z}$. By substituting this in Eq.~\eqref{eq:boa} and then dividing both sides by $N-1$, one finds that
%\begin{eqnarray*}
%\sum_{n=1}^{N-1}{\ln{\left[\Gamma{(x_n)} \cdot \Gamma{(1-x_n)} \right] }} = (N-1) \, \ln{(2 \, \pi)} -\ln{N} \, .
%\label{eq:prodgam}
%\end{eqnarray*}
%By dividing both sides by $N-1$, one has
\begin{eqnarray*}
\sum_{n=1}^{N-1}{\frac{\ln{\Gamma{(x_n)}}+\ln{\Gamma{(1-x_n)}}}{N-1}} = \ln{(2 \, \pi)} -\frac{\ln{N}}{N-1} \, .
%\label{eq:rsum}
\end{eqnarray*}
By taking the limit as $N \rightarrow \infty$, one finds
\begin{equation}
\int_0^{\,1}{\left[\ln{\Gamma{(x)}}+\ln{\Gamma{(1-x)}}\right] \: d x} = \ln{(2 \, \pi)} \, .
\label{eq:intgam}
\end{equation}
Now, note that $\int_0^1{\left[\ln{\Gamma{(x)}}+\ln{\Gamma{(1-x)}}\right] \: d x} = \int_0^1{\ln{\Gamma{(x)}} \: d x} + \int_0^1{\ln{\Gamma{(1-x)}} \: d x}$.  By substituting $y=1-x$ in the latter integral, one finds %that the terms are identical, so
\begin{equation*}
 2 \, \int_0^{\,1}{\ln{\Gamma{(x)}} \: d x} = \ln{(2 \, \pi)} \, , %\int_0^1{\left[\ln{\Gamma{(x)}}+\ln{\Gamma{(1-x)}}\right] \: d x} =
\end{equation*}
which means that
\begin{equation}
\int_0^{\,1}{\ln{\Gamma{(x)}} \: d x} = \frac12 \, \ln{(2 \, \pi)} = \ln{\sqrt{2 \, \pi}} \: .
\end{equation}
This nice evaluation appears on the cover of the gripping book \emph{Irresistible integrals}, by Boros and Moll~\cite{Irresistiveis}.

I left for the reader the less obvious task of using the identity
\begin{equation}
\prod_{n=0}^{N-1} \sin{\left(\pi \frac{n}{N} + \theta \right)} = \frac{\sin(N \theta)}{2^{N-1}} \, ,
\label{eq:teta}
\end{equation}
valid for all positive integer $N$ and all $\theta \in \mathbb{R}$, for evaluating the definite integral
\begin{equation}
\int_0^{\,1}{\ln{\left| \sin{(\pi x + \theta)} \right|} \: d x} \, .
\label{eq:intfim}
\end{equation}
%It may be useful to restrict yourself first to the case $\theta \ne r \, \pi$, with $r \in \mathbb{Q}$.
The proof of the trigonometric identity in Eq.~\eqref{eq:teta} is proposed as an exercise in Ref.~\cite{Niven} (see its Appendix~A.3, Ex.~6).

%\section{Conclusion}
%In this short paper, I have presented a simple method for exactly evaluating the integrals $\int_0^{\pi}{\ln{(\sin{\theta})}\,d\theta}$, $\int_0^{\pi/2}{\ln{(\sin{\theta})}\,d\theta}$, $\int_0^{\pi/2}{\ln{(\cos{\theta})}\,d\theta}$, and $\int_0^{\pi/2}{\ln{(\tan{\theta})}\,d\theta}$ starting from certain finite products of the corresponding trigonometric function at rational multiples of $\pi$. I have shown that the logarithm of each of these products yields sums that can be put in the form of Riemann sums, which yields closed-form results when we take the limit as the number of terms tends to infinity. I have also applied the method for evaluating  $\int_0^1{\ln{\Gamma(x)} \, d x}$ exactly, a task that is otherwise considerably complex.

%\section*{Acknowledgments}
%The author acknowledges a postdoctoral fellowship from CNPq (Brazilian agency).

\end{document}